# EXTREMES OF THE DISCRETE TWO-DIMENSIONAL GAUSSIAN FREE FIELD[1]

By Olivier Daviaud

*Stanford University*

We consider the lattice version of the free field in two dimensions and study the fractal structure of the sets where the field is unusually high (or low). We then extend some of our computations to the case of the free field conditioned on being everywhere nonnegative. For example, we compute the width of the largest downward spike of a given length. Through the prism of these results, we find that the extrema of the free field under entropic repulsion (minus its mean) and those of the unconditioned free field are identical. Finally, when compared to previous results these findings reveal a suggestive analogy between the square of the free field and the two-dimensional simple random walk on the discrete torus.

**1. Introduction.** Let $V_N := \{1, \ldots, N\}^2$, $\partial V_N$ being the points in $V_N$ which have a nearest neighbor outside, and $\mathrm{int}(V_N)$ those which do not. We define the two-dimensional discrete Gaussian free field $\Phi = \{\Phi_x\}_{x \in V_N}$ as follows: $\Phi$ is a family of centered Gaussian random variables with covariance given by the discrete Green's function:

$$(1.1) \qquad G_N(x, y) = \mathbb{E}_x \left( \sum_{i=0}^{\tau_{\partial V_N}} \mathbb{1}_{\eta_i = y} \right), \qquad x, y \in \mathrm{int}(V_N),$$

where $\{\eta_i\}_{i \geq 0}$ is a two-dimensional simple random walk on $\mathbb{Z}^2$, and

$$(1.2) \qquad \tau_{\partial V_N} := \inf\{i \geq 0 : \eta_i \in \partial V_N\}.$$

Received June 2004; revised August 2005.

[1]Supported in part by NSF Grants DMS-04-06042 and DMS-FRG-0244323.

Alternatively, $\Phi$ can be defined as the finite-volume Gibbs measure on $\mathbb{R}^{V_N}$ with Hamiltonian

$$H(\Phi) := \tfrac{1}{8} \sum_{\substack{x,y \in V_N, \\ |x-y|=1}} (\Phi_x - \Phi_y)^2$$

and zero boundary conditions. For any $d \geq 3$ the $d$-dimensional Gaussian free field is defined similarly: the Hamiltonian is the same, and the variance is given by the same formula (using a $d$-dimensional random walk instead).

In statistical physics the Gaussian free field is used to describe the interface between two phases at low temperature. This particular choice for an interface (i.e., a function from $\mathbb{Z}^2$ to $\mathbb{R}$) is motivated in [8], Section 2.4, for example. Using a two-dimensional anisotropic Ising model the author explains that when one of the parameters becomes large, the measure concentrates on configurations where particles with the same spin form a connected component, and where the separation line between the two components has minimal horizontal length (and is therefore a function). Another explanation can be found in [7], Section 1.3, where it is pointed out that the Hamiltonian for the two-dimensional Ising model is proportional to the perimeter of the interface. At low temperature this leads again to the types of configurations described above. In addition, part of the Gaussian free field's appeal is explained by its tractability. Indeed, it is the only massless model with nearest-neighbor Hamiltonian in dimension $d \geq 2$ (see [8], Sections 2.1, 2.2 and 2.3 for definitions) for which the distribution is known [a centered Gaussian with covariance given by (1.1)]. Yet the model is rich enough to illustrate many interesting phenomena. In this paper we will focus on two of them: the multiscale structure of the field, and the phenomenon of entropic repulsion.

The term *entropic repulsion* refers to the behavior of the field when it is conditioned on being everywhere nonnegative (i.e., when the interface is in presence of a hard wall). In that case, it is known that the field is pushed very far away as $N$ becomes large (see [2] for the case $d=2$, and [3] for $d \geq 3$), while at the same time the deformation undergone by the field is very small. To quote a sentence from Giacomin in his survey [8], "the free field does not give up easily its freedom of fluctuating." In dimension $d \geq 3$, these observations are quantified in [3] where a precise study of the Gaussian free field under entropic repulsion is undertaken. It is proved in that paper that the conditioned free field converges weakly toward the original field shifted by a constant (which varies with $N$). However, it is easily seen in the context of [6] that the extrema of the shifted and conditioned field behave differently (their minima are different, e.g.). In this paper we argue that in two dimensions this is not the case: the extrema of the conditioned and shifted free field look alike. We will use this "guess" as a strategy and deduce results on the conditioned free field from results on the original free field. This leads us to our first theorem which gives the typical width of a downward spike of a given (negative) height for the conditioned free field. By "downward spike" we mean the conditioned free field restricted to a set where it is uniformly unusually close to 0. Once and for all we fix $l \in (0, 1/2)$, and define $V_N^l := \{x \in V_N : \text{dist}(x, V_N^c) \geq lN\}$, $\Omega_{N,l}^+ := \{\Phi_x \geq 0, x \in V_N^l\}$.

By *conditioned free field* (or CFF) we denote the law of $\Phi$ conditioned on $\Omega_{N,l}^+$. Finally, in what follows $B(x,a)$ denotes the box of center $x$ and edge length $a$, and $g := 2/\pi$. Then:

THEOREM 1.1. *Let $0 < \eta < 1$, $l > 0$ and*

$$(1.3) \quad D_N^+(\eta) := \sup\left\{a \in \mathbb{N} \ such \ that \ \exists x \in V_N^l : \max_{B(x,a)} \Phi \leq 2(1-\eta)\sqrt{g}\log N\right\}.$$

*Then*

$$(1.4) \quad \lim_{N \to \infty} \frac{\log D_N^+(\eta)}{\log N} = \frac{1}{2} - \frac{\eta}{2} \qquad in \ probability \ under \ \mathbf{P}(\cdot|\Omega_{N,l}^+).$$

Said in words, since it is proved in [2] that the mean of the CFF is asymptotically equal to $2\sqrt{g}\log N$, Theorem 1.1 states that the largest downward spike of height $2\sqrt{g}\eta\log N$ has width $N^{1/2-\eta/2}$. The total number of points that are at or below this low level is given by the following.

PROPOSITION 1.2. *Let $\mathcal{L}_N^+(\eta) := \{x \in V_N^l : \Phi_x \leq 2\sqrt{g}(1-\eta)\log N\}$. Then for $0 < \eta < 1$,*

$$(1.5) \quad \lim_{N \to \infty} \frac{\log |\mathcal{L}_N^+(\eta)|}{\log N} = 2(1-\eta^2) \qquad in \ probability \ under \ \mathbf{P}(\cdot|\Omega_{N,l}^+).$$

To obtain results such as Theorem 1.1, one first needs to understand the behavior of the extrema of the free field without conditioning. The starting point in this direction is [2] where the maximum of the Gaussian free field was computed for the first time. It is proved in that paper that

$$(1.6) \quad \lim_{N \to \infty} \mathbf{P}\left(\left|\sup_{x \in V_N^l} \Phi_x - 2\sqrt{g}\log N\right| \geq \delta\log N\right) = 0.$$

The proof relies on the idea that the easiest way for the field to attain its maximum is to have an upward shift on all scales $N^\alpha$, $0 < \alpha < 1$. We will exploit this *multiscale decomposition* to study the repartition of the points where the field is unusually high. As we will see below, the set of such points exhibits a fractal structure; this structure is specific to the two-dimensional free field, as it vanishes in lower or higher dimension (i.e., the exponents in our theorems—which are discrete analogues of the Hausdorf dimension—become trivial when $d \neq 2$). Maybe more importantly, this structure is identical to the one uncovered in [5], by which many of our results—and many of our proofs—are inspired. More precisely, if $x$ is a point on the two-dimensional discrete torus of size $N$ [i.e., $(\mathbb{Z}/N\mathbb{Z})^2$], and if $\tau(x)$ denotes the first time $x$ is visited by a simple random walk starting at 0, then our next five theorems are shown in [5] to hold when replacing $\Phi$ with $\sqrt{2\tau}/N$. It is an open problem to further explore this correspondence.

We now need to specify what we mean by "unusually high." We say that $x \in V_N$ is an *$\alpha$-high point* for the GFF if $\Phi_x \geq 2\sqrt{g}\alpha \log N$. We start by computing the number of such high points:

THEOREM 1.3. *Let* $\mathcal{H}_N(\eta) := \{x \in V_N^l : \Phi_x \geq 2\sqrt{g}\eta \log N\}$ *be the set of $\eta$-high points. Then for $0 < \eta < 1$,*

$$(1.7) \qquad \lim_{N \to \infty} \frac{\log |\mathcal{H}_N(\eta)|}{\log N} = 2(1 - \eta^2) \qquad \text{in probability.}$$

*Moreover, for all $\delta > 0$ there exists a constant $c > 0$ such that*

$$(1.8) \qquad \mathbf{P}(\{|\mathcal{H}_N(\eta)| \leq N^{2(1-\eta^2)-\delta}\}) \leq \exp[-c(\log N)^2]$$

*for $N$ large enough.*

These points are not evenly spread, and typically appear in cluster. This is quantified by the two following theorems which give the number of high points in a small fixed neighborhood and around a high point, respectively.

THEOREM 1.4. *Let* $D(x, \rho) := \{y \in V_N : |y - x| \leq \rho\}$ *be the disk of radius $\rho$ centered at $x$. For $0 < \alpha < \beta < 1$, and $\delta > 0$,*

$$(1.9) \quad \lim_{N \to \infty} \max_{x \in V_N^l} \mathbf{P}\left(\left|\frac{\log |\mathcal{H}_N(\alpha) \cap D(x, N^\beta)|}{\log N} - 2\beta(1 - (\alpha/\beta)^2)\right| > \delta\right) = 0.$$

THEOREM 1.5. *For $0 < \alpha$, $\beta < 1$ and $\delta > 0$,*

$$(1.10) \quad \begin{aligned} \lim_{N \to \infty} \max_{x \in V_N^l} \mathbf{P}\bigg(&\left|\frac{\log |\mathcal{H}_N(\alpha) \cap D(x, N^\beta)|}{\log N}\right. \\ & \left. - 2\beta(1 - \alpha^2)\right| > \delta \,\Big|\, x \in \mathcal{H}_N(\alpha)\bigg) = 0. \end{aligned}$$

Note that the exponents we obtain are, respectively, lower and higher than if the high points were evenly spread (indeed in that case, in view of Theorem 1.3, there should be $N^{2\beta-2\alpha^2}$ $\alpha$-high points in a neighborhood of size $N^\beta$). An even more complex structure is obtained when one considers the pairs of high points. Their number is given by:

THEOREM 1.6. *Let $0 < \alpha, \beta < 1$. Then*

$$(1.11) \qquad \lim_{N \to \infty} \frac{\log |\{(x, y) \in \mathcal{H}_N(\alpha) : d(x, y) \leq N^\beta\}|}{\log N} = \rho(\alpha, \beta)$$

*in probability, where*

$$(1.12) \qquad \rho(\alpha, \beta) = 2 + 2\beta - 2\alpha^2 \inf_{\gamma \in \Gamma_{\alpha, \beta}} F_{2, \beta}(\gamma),$$

$$(1.13) \qquad F_{h, \beta}(\gamma) = \gamma^2 (1 - \beta) + \frac{h(1 - \gamma(1 - \beta))^2}{\beta},$$

$$(1.14) \qquad \Gamma_{\alpha, \beta} = \{\gamma \geq 0 : 2 - 2\beta - 2\alpha^2 F_{0, \beta}(\gamma) \geq 0\}.$$

This result should be compared to the mean numbers of pairs. This last quantity is easy to compute using Lemma 8.2 and is of the order of $N^{\tilde\rho(\alpha, \beta)}$ where

$$\tilde\rho(\alpha, \beta) = 2 + 2\beta - 2\alpha^2 \inf_\gamma F_{2, \beta}(\gamma).$$

Therefore the mean and median numbers of pairs have different orders of magnitude. Finally, in order to prove Theorem 1.1 we will need a similar result for the unconditioned free field. The following theorem gives the width of the largest upward spike of a given length:

THEOREM 1.7. *Let $-1 < \eta < 1$, $l > 0$ and let $D_N(\eta)$ be the side length of the biggest square where the GFF is uniformly greater than $2\eta\sqrt{g}\log N$, that is,*

$$(1.15) D_N(\eta) := \sup \left\{ a \in \mathbb{N} \text{ such that } \exists x \in V_N^l : \min_{B(x, a)} \Phi \geq 2\eta\sqrt{g}\log N \right\}.$$

*Then*

$$(1.16) \qquad \lim_{N \to \infty} \frac{\log D_N(\eta)}{\log N} = \frac{1}{2} - \frac{\eta}{2} \qquad \text{in probability.}$$

We will first prove each result on the (unconditioned) GFF in a different section, in the order in which they were introduced. Section 7 contains the proofs of Theorem 1.1 and Proposition 1.2. Finally, the last section gathers the technical lemmas that will be needed along the way.

REMARKS. (i) Since the two-dimensional free field, once properly scaled, converges to a continuous version (see [12] or Section 9 for an introduction to this object), one could ask whether our results carry over to that case. In [9], the concept of high point is defined for the continuous free field, and an analogue of Theorem 1.3 is proved (i.e., the Hausdorff dimension of the set of high points is computed).

(ii) Similarly, it is proved in [10] that the height function associated to a random domino tiling of $V_N$, for example, converges to the continuous version of the Gaussian free field. It would therefore be interesting to know whether our results hold in this context as well.

**2. Theorem 1.3: Number of high points.** We start with the following lemma, which gives some estimates on the covariance structure of the GFF. These estimates will be used extensively in the sequel.

LEMMA 2.1. *There exist $c$ and $c(l)$ (independent of $N$) such that*

$$(2.1) \qquad \sup_{x \in V_N^l} \operatorname{Var}(\Phi_x) \le g \log N + c,$$

$$(2.2) \qquad \sup_{\substack{x,y \in V_N^l \\ x \ne y}} \left( \operatorname{Cov}(\Phi_x, \Phi_y) - g(\log N - \log |y - x|) \right) \le c$$

*and*

$$(2.3) \qquad \sup_{x \in V_N^l} |\operatorname{Var}(\Phi_x) - g \log N| \le c(l),$$

$$(2.4) \qquad \sup_{\substack{x,y \in V_N^l \\ x \ne y}} |\operatorname{Cov}(\Phi_x, \Phi_y) - g(\log N - \log |y - x|)| \le c(l).$$

PROOF. These results are easy consequences of [11], Theorem 1.6.6, Proposition 1.6.7. To see this, let $C_N := D(0, N)$ be the disk of center $0$ and radius $N$, and let $\tau_{\partial C_N}$ denote the first exit time of such a disk by a two-dimensional simple random walk $\{\eta_i, \ i \ge 0\}$ started at $0$. With this notation we define

$$G_{C_N}(x, y) := \mathbb{E}_x \left( \sum_{i=0}^{\tau_{\partial C_N}} \mathbb{1}_{\eta_i = y} \right), \qquad x, y \in \operatorname{int}(V_N).$$

Then, for $x \in V_N$,

$$\begin{aligned}
\operatorname{Var}(\Phi_x) &= G_N(x, x) \\
&\le G_{C_{2N}}(0, 0) \\
&= g \ln 2N + k + O(1/(2N))
\end{aligned}$$

where the first inequality follows from $D(x, 2N) \supset V_N$, and the last line is [11], Theorem 1.6.6 ($k$ is a constant). Therefore $G_N(x, x) \le g \log N + c$ for some universal constant $c$. This proves (2.1). To derive (2.3) we note that for $x \in V_N^l$,

$$\begin{aligned}
\operatorname{Var}(\Phi_x) &= G_N(x,x) \\
&\geq G_{C_{lN}}(0,0) \\
&= g \ln lN + k + O(1/(lN)) \\
&\geq g \ln N + c(l)
\end{aligned}$$

for some constant $c(l)$ depending only on $l$. Equations (2.2) and (2.4) are derived from [11], Proposition 1.6.7, in a similar fashion. $\square$

PROOF OF THE UPPER BOUND. This part of the theorem simply follows from Chebyshev's inequality:

$$\begin{aligned}
\mathbf{P}(|\mathcal{H}_N(\eta)| \geq N^{2(1-\eta^2)+\delta}) &\leq N^{-2(1-\eta^2)-\delta} \mathbb{E}(|\mathcal{H}_N(\eta)|) \\
&\leq N^{-2(1-\eta^2)-\delta} N^2 \max_{x \in V_N} \mathbf{P}(\Phi_x \geq 2\eta\sqrt{g}\log N) \\
&\leq N^{-2(1-\eta^2)-\delta} N^2 \exp\left[\frac{-4\eta^2 g(\log N)^2}{2g\log N + c}\right] \\
&\leq N^{-2(1-\eta^2)-\delta/2} N^{2-2\eta^2} \longrightarrow 0
\end{aligned}$$
(2.5)

where (2.5) is a consequence of (2.1) combined with the following well-known tail estimate: if $X$ is a standard normal variable, then for any $a \geq 0$,

$$\mathbf{P}(|X| \geq a) \leq e^{-a^2/2}. \tag{2.6}$$

$\square$

PROOF OF THE LOWER BOUND. We start with two remarks that will remain in effect throughout this paper. First, for ease of exposition we ignore all discretization issues. An illustrating example is that by "box $B$ of center $x \in V_N$ and side length $a$" ["$B(x,a)$"] we mean a box of center $x$ and side length $b$, where $b$ is an even integer close to $a$, and is adjusted depending on the context. For example, $b$ can be chosen so as to make $B$ adjacent to some other box. Second, we will make extensive use of the notation introduced in [2]. We recall it now for the reader's convenience. For any subset $A \subseteq V_N$, $\mathcal{F}_A$ will denote the $\sigma$-algebra generated by $\Phi_x$, $x \in A$. For any box $B$, we write $x_B$ for the center of $B$ and let $\Phi_B := \mathbb{E}(\Phi_{x_B} | \mathcal{F}_{\partial B})$. For $0 \leq \alpha < 1$, we let $\Pi_\alpha$ be the collection of adjacent sub-boxes of edge length $N^\alpha$ in $V_N^l$ such that $V_N^l = \bigcup_{B \in \Pi_\alpha} B$. In particular, $\Pi_0$ is defined as $\{\{x\} : x \in V_N^l\}$, and for $\{x\} \in \Pi_0$, $\Phi_{\{x\}}$ is simply $\Phi_x$. Next, we define a collection of boxes as follows. We fix $1/2 < \alpha < 1$, $K$ an integer greater than 2, and let $\alpha_i := \alpha(K-i+1)/K$, $1 \le i \le K+1$. We first let $\Gamma_{\alpha_1} := \Pi_{\alpha_1}$. Then assuming that $\Gamma_{\alpha_i}$ has been defined, for any $B \in \Gamma_{\alpha_i}$ we draw a square of side length $N^{\alpha_i}/2$ with the same center as $B$. The collection of sub-boxes in $\Pi_{\alpha_{i+1}}$ that intersect that square is called $\Gamma_{B,\alpha_{i+1}}$. Note that $|\Gamma_{B,\alpha_{i+1}}| = N^{2\alpha/K}/4$ by volume considerations. Finally, we let $\Gamma_{\alpha_{i+1}} := \bigcup_{B \in \Gamma_{\alpha_i}} \Gamma_{B,\alpha_{i+1}}$.

The lower bound of (1.7) as well as (1.8) follow from a simple modification of the proof of [2], Theorem 2b, which we now explain. First let

$$D_1 := \{B \in \Pi_{\alpha_1} : \Phi_B \ge 0\}$$

and let $C_1$ be the event $\{\#D_1 \ge N^\kappa\}$, where we choose $\kappa > 0$ such that

$$(2.7) \qquad \mathbf{P}(C_1) \ge 1 - \exp(-a(\log N)^2)$$

for some $a > 0$. That such a $\kappa$ exists is guaranteed by [2], Lemma 8. Next for $2 \le k \le K+1$ we let $n_k := N^{\kappa + (2\alpha(k-1)/K)(1-\eta^2)}$ and, after setting $\gamma_K := 1/K$, we define

$$(2.8) \quad \begin{aligned} D_k &:= \{\underline{B}^{(k)} : \Phi_{B_i} \ge (\alpha - \alpha_i)\eta 2\sqrt{g}(1-\gamma_K)\log N, \ \forall 1 \le i \le k\}, \\ C_k &:= \{\#D_k \ge n_k\}, \end{aligned}$$

where $\underline{B}^{(k)}$ denotes a sequence of boxes $(B_1, \ldots, B_k)$ satisfying $B_1 \supset B_2 \supset \cdots \supset B_k$, and $B_i \in \Gamma_{\alpha_i}$. On the event $C_k$ there are at least $n_k$ sequences $\underline{B}^{(k)} = (B_{j,1}, B_{j,2}, \ldots, B_{j,k})$, $1 \le j \le n_k$, in the set $D_k$. Then for $k \le K$

$$C_k \cap C_{k+1}^c \subseteq C_k \cap \left\{ \sum_{j=1}^{n_k} \zeta_j \le \frac{4n_{k+1}}{N^{2\alpha/K}} \right\},$$

where

$$\zeta_j := \frac{1}{|\Gamma_{B_{j,k},\alpha_{k+1}}|} \sum_{B \in \Gamma_{B_{j,k},\alpha_{k+1}}} \mathbb{1}_{\{\Phi_B - \Phi_{B_{j,k}} \ge \eta\alpha(2\sqrt{g}/K)(1-\gamma_K)\log N\}}.$$

Let

$$A_k = \bigcup_{A \in \Gamma_{\alpha_k}} \bigcup_{B \in \Gamma_{A,\alpha_{k+1}}} \{|\Phi_A - \mathbb{E}(\Phi_B|\mathcal{F}_{\alpha_k})| \ge \varepsilon \eta\alpha(2\sqrt{g}/K)(1-\gamma_K)\log N\},$$

where $\varepsilon > 0$, and for any $0 < \beta < 1$ we let $\mathcal{F}_\beta := \sigma(\Phi_x : x \in \bigcup_{B \in \Pi_\beta} \partial B)$. Since by [2], Lemma 12,

$$(2.9) \qquad \mathrm{Var}(\Phi_A - \mathbb{E}(\Phi_B|\mathcal{F}_{\alpha_k})) \le c$$

for some universal constant $c$ and any $A \in \Gamma_{\alpha_k}$, $B \in \Gamma_{A,\alpha_{k+1}}$, the union bound and (2.6) give

$$(2.10) \qquad \mathbf{P}(A_k) \le e^{-d(\log N)^2}$$

for some $d > 0$. Note that on $A_k^c$, $\zeta_j \geq \tilde{\zeta}_j$ where

$$\tilde{\zeta}_j = \frac{1}{|\Gamma_{B_{j,k},\alpha_{k+1}}|} \sum_{B \in \Gamma_{B_{j,k},\alpha_{k+1}}} \mathbb{1}_{\{\Phi_B - \mathbb{E}(\Phi_B|\mathcal{F}_{\alpha_k}) \geq \eta\alpha(2\sqrt{g}/K)(1-\gamma_K)(1+\varepsilon)\log N\}}.$$

Let $\mathrm{Var}_{\mathcal{F}_{\alpha_k}}$ denote the conditional variance given $\mathcal{F}_{\alpha_k}$. If $B \in \Gamma_{B_{j,k},\alpha_{k+1}}$ we have

$$\begin{align}
(2.11) \qquad \mathrm{Var}_{\mathcal{F}_{\alpha_k}}(\Phi_B) &= \mathrm{Var}_{\mathcal{F}_{\alpha_k}}(\mathbb{E}(\Phi_{x_B}|\mathcal{F}_{\partial B})) \\
&= \mathrm{Var}_{\mathcal{F}_{\alpha_k}}(\Phi_{x_B}) - \mathbb{E}(\mathrm{Var}_{\mathcal{F}_{\partial B}}(\Phi_{x_B})|\mathcal{F}_{\alpha_k}).
\end{align}$$

It can be seen that under $\mathbf{P}(\cdot|\mathcal{F}_{\partial B})$, the distribution of $\{\Phi_x\}_{x \in B}$ is the free field on $B$ with boundary conditions given by $\Phi_{|\partial B}$. That is, under $\mathbf{P}(\cdot|\mathcal{F}_{\partial B})$, $\{\Phi_x\}_{x \in B}$ is a Gaussian variable with covariance given by $G_{N^{\alpha_{k+1}}}(\cdot, \cdot)$ [with $B$ replaced with $V_N$ in (1.1)] and mean a discrete harmonic function with boundary conditions $\{\Phi_x\}_{x \in \partial B}$ (cf. [2], Section 2.3 for more details). In particular, $\mathrm{Var}_{\mathcal{F}_{\partial B}}(\Phi_{x_B})$ is constant and can be bounded using Lemma 2.1. These considerations also apply to $\mathrm{Var}_{\mathcal{F}_{\alpha_k}}(\Phi_{x_B})$, and when combined with (2.11) give

$$\begin{align}
(2.12) \qquad \mathrm{Var}_{\mathcal{F}_{\alpha_k}}(\Phi_B) &\geq \alpha_k g \log N - \alpha_{k+1} g \log N - 2c \\
&= \frac{g\alpha}{K} \log N - 2c.
\end{align}$$

Now the random variables $\tilde{\zeta}_j$ are i.i.d. (under $\mathbf{P}$), and by the discussion above

$$(2.13) \qquad \mathbb{E}(\tilde{\zeta}_j|\mathcal{F}_{\alpha_k}) \geq \mathbf{P}(X \geq \eta\alpha(2\sqrt{g}/K)(1-\gamma_K)(1+\varepsilon)\log N)$$

where $X$ is centered Gaussian with variance greater than

$$(2.14) \qquad \frac{g\alpha}{K} \log N - 2c$$

for some $c > 0$. Next recall that if $Y$ is a standard normal variable, for any $a \geq 1$

$$(2.15) \qquad \mathbf{P}(|Y| \geq a) \geq \frac{e^{-a^2/2}}{\sqrt{2\pi}a}.$$

This fact, together with (2.13), (2.14) and $1 > 1 - \gamma_K > 0$, gives

$$\mathbb{E}(\tilde{\zeta}_j|\mathcal{F}_{\alpha_k}) \geq N^{-(2/K)(1-\gamma_K)\eta^2\alpha(1+\varepsilon)^2}$$

provided $N$ is large enough. Consequently on the event $C_k \cap A_k^c$

$$\begin{align}
C_{k+1}^c &\subseteq \left\{ \sum_{j=1}^{n_k}(\tilde{\zeta}_j - \mathbb{E}(\tilde{\zeta}_j|\mathcal{F}_{\alpha_k})) \leq \frac{4n_{k+1}}{N^{2\alpha/K}} - n_k N^{-(2/K)(1-\gamma_K)\eta^2\alpha(1+\varepsilon)^2} \right\} \\
&\subseteq \left\{ \left| \sum_{j=1}^{n_k}(\tilde{\zeta}_j - \mathbb{E}(\tilde{\zeta}_j|\mathcal{F}_{\alpha_k})) \right| \geq \frac{1}{2}n_k N^{-(2/K)(1-\gamma_K)\eta^2(1+\varepsilon)^2\alpha} \right\}
\end{align}$$

provided

$$-\frac{2}{K}(1-\gamma_K)\eta^2\alpha(1+\varepsilon)^2 > -\frac{2\alpha}{K}\eta^2,$$

that is,

$$(1+\varepsilon)^2(1-\gamma_K) < 1.$$

We then bound the probability of the above set using [2], Lemma 11:

LEMMA. *Let* $Z_1, \ldots, Z_n$ *be i.i.d. real-valued random variables satisfying* $\mathbb{E}(Z_i) = 0$, $\sigma^2 = \mathbb{E}Z_i^2$, $\|Z_i\|_\infty \le 1$. *Then for any* $t > 0$

$$\mathbf{P}\left(\left|\sum_{i=1}^n Z_i\right| \ge t\right) \le 2\exp\left[-\frac{t^2}{2n\sigma^2 + 2t/3}\right].$$

By taking $n = n_k$, $t = \frac{1}{2}nN^{-(2/K)(1-\gamma_K)\eta^2(1+\varepsilon)^2\alpha}$ and $\sigma \le 1$ we obtain that on $A_k^c \cap C_k$

$$(2.16) \quad \mathbf{P}(C_{k+1}^c | \mathcal{F}_{\alpha_k}) \le 2\exp[-cN^{\kappa + (2\alpha/K)(1-\eta^2)(k-1) - (4/K)(1-\gamma_K)\eta^2(1+\varepsilon)^2\alpha}].$$

Now

$$\mathbf{P}(C_{K+1}^c) \le \sum_{k=2}^{K+1}(\mathbf{P}(C_k^c \cap C_{k-1} \cap A_{k-1}^c) + \mathbf{P}(A_{k-1})) + \mathbf{P}(C_1^c)$$

which together with (2.7), (2.10), (2.16) and $\varepsilon$ small enough implies that

$$\mathbf{P}(C_{K+1}^c) \le \exp(-c(\log N)^2)$$

for some constant $c > 0$. Since

$$C_{K+1} \subset \{|\mathcal{H}_N(\eta\alpha(1-\gamma_K))| \ge n_{K+1} := N^{\kappa + 2\alpha(1-\eta^2)}\},$$

taking $\alpha$ close to 1 and $K$ large enough completes the proof of the lower bound. $\square$

## 3. Theorem 1.4: High points in a small neighborhood.

PROOF OF THEOREM 1.4. Let $B := B(x, 4N^\beta)$. The idea behind the proof (inspired by [5]) is that typically $\Phi_B$ is close to 0 (i.e., with the notation below: $\mathbf{P}(D_+ \cap D_-)$ is close to 1). Therefore, since for $y \in D(x, N^\beta)$, $\mathbb{E}(\Phi_y | \mathcal{F}_{\partial B})$ is close to $\Phi_B$ [i.e., in the notation below: $\mathbf{P}(A)$ is close to 0], for such $y$ we have that

$$(3.1) \qquad \Phi_y \approx \Phi_y - \mathbb{E}(\Phi_y | \mathcal{F}_{\partial B}).$$

Since $\alpha \log N \approx (\alpha/\beta) \log 4N^\beta$, in view of (3.1) the event $\{y \in \mathcal{H}_N(\alpha)\}$ roughly corresponds to $\{\Phi_y - \mathbb{E}(\Phi_y | \mathcal{F}_{\partial B}) \geq (\alpha/\beta) \log 4N^\beta\}$. We conclude by conditioning on $\mathcal{F}_{\partial B}$ and by using Theorem 1.3 which gives the number of such points $y$. Indeed, as noticed in the previous section, under $\mathbf{P}(\cdot | \mathcal{F}_{\partial B})$, the family $\{\Phi_y - \mathbb{E}(\Phi_y | \mathcal{F}_{\partial B}), \ y \in B\}$ has the same law as the free field on $V_{4N^\beta}$.

Fix $\eta, \delta > 0$ and let

$$D_+ := \{\Phi_B \leq \eta 2\sqrt{g} \log N\},$$

$$C_+ := \{|\mathcal{H}_N(\alpha) \cap D(x, N^\beta)| \geq N^{2\beta(1-(\alpha/\beta)^2)+\delta}\}$$

and

$$A := \bigcup_{y \in D(x, N^\beta)} \{|\mathbb{E}(\Phi_y | \mathcal{F}_{\partial B}) - \Phi_B| \geq \varepsilon 2\sqrt{g} \log N\}$$

where $\varepsilon > 0$. By [2], Lemma 12, $\mathbf{P}(A)$ tends to 0. So does $\mathbf{P}(D_+^c)$ [by (2.6) and Lemma 2.1] and therefore

$$\begin{aligned}
(3.2) \qquad \mathbf{P}(C_+) &= \mathbb{E}(\mathbf{P}(C_+ | \mathcal{F}_{\partial B})) \\
&\leq \mathbf{P}(A) + \mathbf{P}(D_+^c) + \mathbb{E}(\mathbf{P}(C_+ | \mathcal{F}_{\partial B}) \mathbb{1}_{A^c} \mathbb{1}_{D_+}) \\
&\leq o(1) + \mathbf{P}\left(\left|\mathcal{H}_{4N^\beta}\left(\frac{\alpha - \varepsilon'}{\beta}\right)\right| \geq N^{2\beta(1-(\alpha/\beta)^2)+\delta}\right)
\end{aligned}$$

where $\varepsilon'$ is such that

$$\frac{\alpha - \varepsilon'}{\beta} \log 4N^\beta = (\alpha - \varepsilon - \eta) \log N.$$

If $N, \eta$ and $\varepsilon$ are, respectively, large, small and small enough so that

$$2\beta\left(1 - \left(\frac{\alpha - \varepsilon'}{\beta}\right)^2\right) < 2\beta\left(1 - \left(\frac{\alpha}{\beta}\right)^2\right) + \delta,$$

then it follows from Theorem 1.3 that the right-hand side of (3.2) tends to 0, which proves the upper bound.

We prove the lower bound in a similar way: we first define

$$D_- := \{\Phi_B \geq -\eta 2\sqrt{g} \log N\},$$

$$C_- := \{|\mathcal{H}_N(\alpha) \cap D(x, N^\beta)| \leq N^{2\beta(1-(\alpha/\beta)^2)-\delta}\},$$

and introduce the notation

$$(3.3) \qquad \mathcal{H}_N^s(\alpha) := \{x \in V_N^s : \Phi_x \geq 2\sqrt{g}\eta \log N\},$$

where $0 < s < 1/2$. As before

$$\mathbf{P}(C_-) = \mathbb{E}(\mathbf{P}(C_- | \mathcal{F}_{\partial B}))$$

$$\leq \mathbf{P}(A) + \mathbf{P}(D_-^c) + \mathbb{E}(\mathbf{P}(C_-|\mathcal{F}_{\partial B})\mathbb{1}_{A^c}\mathbb{1}_{D_-})$$

$$= o(1) + \mathbb{E}(\mathbf{P}(C_-|\mathcal{F}_{\partial B})\mathbb{1}_{A^c}\mathbb{1}_{D_-})$$

$$\leq o(1) + \mathbf{P}\left(\left|\mathcal{H}_{4N^\beta}^{3/8}\left(\frac{\alpha+\varepsilon}{\beta}\right)\right| \leq N^{2\beta(1-(\alpha/\beta)^2)-\delta}\right)$$

where $\varepsilon'$ is chosen so that

$$\frac{\alpha+\varepsilon'}{\beta}\log 4N^\beta = (\alpha+\varepsilon+\eta)\log N.$$

We then conclude as in the upper bound case. $\quad\square$

## 4. Theorem 1.5: Clusters of high points.

PROOF OF THE LOWER BOUND. The proof (inspired by [5]) roughly goes as follows. As in the previous section, we pick a box $B$ of size $4N^\beta$ centered at $x$. This time, since we condition on $\{x \in \mathcal{H}_N(\alpha)\}$, the typical value of $\Phi_B$ is $(1-\beta)2\alpha\sqrt{g}\log N$ (instead of $0$). Since again $\mathbb{E}(\Phi_y|\mathcal{F}_{\partial B}) \approx \Phi_B$ for $y$ in $D(x, N^\beta)$, $\{y \in \mathcal{H}_N(\alpha)\}$ if and only if $\{\Phi_y - \mathbb{E}(\Phi_y|\mathcal{F}_{\partial B}) \geq \alpha\beta 2\sqrt{g}\log N\}$. Applying Theorem 1.3 after conditioning on $\mathcal{F}_{\partial B}$ gives the number of such points $y$.

We start with a few definitions. We let $\eta, \delta > 0$, $B := B(x, 4N^\beta)$,

$$b^-(\alpha, \beta, \eta, N) := 2\sqrt{g}(\alpha(1-\beta) - \eta)\log N,$$

$$E := \{|\mathcal{H}_N(\alpha) \cap D(x, N^\beta)| \leq N^{2\beta(1-\alpha^2)-\delta}\},$$

$$F := \{\Phi_B \geq b^-(\alpha, \beta, \eta, N)\},$$

and finally $G := \{x \in \mathcal{H}_N(\alpha)\}$. Then

$$\mathbf{P}(E|G) = o(1) + (1+o(1))\mathbf{P}(E|F, G)$$

follows from Lemma 8.1. Next, by Cauchy–Schwarz's inequality,

$$\mathbf{P}(E|F, G) = \frac{\mathbf{P}(E, F, G)}{\mathbf{P}(F, G)} \leq \frac{1}{\mathbf{P}(F, G)}\sqrt{\mathbf{P}(G)\mathbf{P}(E, F)}$$

$$= \frac{1}{\mathbf{P}(F|G)\mathbf{P}(G)}\sqrt{\mathbf{P}(G)\mathbf{P}(F)\mathbf{P}(E|F)}$$

$$= (1+o(1))\sqrt{\frac{\mathbf{P}(F)}{\mathbf{P}(G)}\mathbf{P}(E|F)}$$

where the last step follows from Lemma 8.1. Since we know by Lemma 2.1 and (2.15) that $\mathbf{P}(F)/\mathbf{P}(G) \leq \exp(d\log N)$ for some $d$, the lower bound will follow once we prove that $\mathbf{P}(E|F) \leq \exp(-c(\log N)^2)$ for some $c > 0$. Now let

$$A := \bigcup_{y \in B(x, N^\beta)} \{|\mathbb{E}(\Phi_y|\mathcal{F}_{\partial B}) - \Phi_B| \geq 2\sqrt{g}\varepsilon \log N\}.$$

Equation (2.15), combined with a rerun of the proof of (2.12), yields that $\mathbf{P}(F) \geq \exp[-d\log N]$ for some $d > 0$. It also follows from [2], Lemma 12 that $\mathbf{P}(A) \leq \exp[-c(\log N)^2]$ for some $c > 0$. Hence $\mathbf{P}(A|F)$ is lower than $\exp[-c(\log N)^2]$ for some other constant $c > 0$. Thus, $\mathbf{P}(E|F)$ is lower than

$$\exp[-c(\log N)^2] + \frac{\mathbb{E}(\mathbf{P}(E|\mathcal{F}_{\partial B})\mathbb{1}_{A^c}\mathbb{1}_F)}{\mathbf{P}(F)}.$$

Now on $A^c \cap F$,

$$
\begin{aligned}
(4.1) \qquad \mathbf{P}(|\mathcal{H}_N(\alpha) \cap D(x, N^\beta)| &\leq N^{2\beta(1-\alpha^2)-\delta}|\mathcal{F}_{\partial B}) \\
&\leq \mathbf{P}(|\mathcal{H}^{3/8}_{4N^\beta}(\alpha + \varepsilon')| \leq N^{2\beta(1-\alpha^2)-\delta})
\end{aligned}
$$

where $\varepsilon'$ is such that

$$(\alpha - (\alpha(1-\beta) - \eta) + \varepsilon)\log N = (\alpha + \varepsilon')\log 4N^\beta.$$

Now from (1.8) we know (4.1) is lower than $\exp[-c(\log N)^2]$ for some $c > 0$ provided $\varepsilon'$ is small enough (which is achieved by taking $\eta$, $\varepsilon$ and $N$ small, small and large enough, resp.). This concludes the proof of the lower bound. □

PROOF OF THE UPPER BOUND. For the upper bound, we condition on $\Phi_B$ and use a first moment method. We let $K$ be a positive integer and $\beta_j := \frac{j}{K}\beta$ for $j = 1, \ldots, K$. We also define $D_1 := D(x, N^{\beta_1})$ and

$$D_j := D(x, N^{\beta_j}) - D(x, N^{\beta_{j-1}}), \qquad j = 1, \ldots, K.$$

We have $D(x, N^\beta) = \bigcup_{j=1}^{K} D_j$. Then

$$\{|\mathcal{H}_N(\alpha) \cap D(x, N^\beta)| \geq N^{2\beta(1-\alpha^2)+\varepsilon}\}$$

$$\subseteq \bigcup_{j=1}^{K}\{|\mathcal{H}_N(\alpha) \cap D_j| \geq N^{2\beta_j(1-\alpha^2)+\varepsilon/2}\}$$

for $N$ large enough. So we just need to prove that

$$\mathbf{P}(|\mathcal{H}_N(\alpha) \cap D_j| \geq N^{2\beta_j(1-\alpha^2)+\varepsilon/2}|x \in \mathcal{H}_N(\alpha)) \xrightarrow{N\to\infty} 0.$$

Note that we can restrict ourselves to values of $\beta_j$ such that

$$2\beta_j(1-\alpha^2) + \varepsilon/2 \leq 2\beta_j,$$

that is, $\beta_j \geq \frac{\varepsilon}{4\alpha^2}$ $(> 0)$. Let $B_j := B(x, 4N^{\beta_j})$ and

$$C := \{|\mathcal{H}_N(\alpha) \cap D_j| \geq N^{2\beta_j(1-\alpha^2)+\varepsilon/2}\}.$$

Finally let $b^+(\alpha, \beta_j, \eta, N) := 2\sqrt{g}\alpha(1 - \beta_j + \eta) \log N$. Then, using Lemma 8.1,

$$\mathbf{P}(C|x \in \mathcal{H}_N(\alpha)) = \mathbf{P}(C \cap \{\Phi_{B_j} \leq b^+(\alpha, \beta, \eta, N)\}|x \in \mathcal{H}_N(\alpha)) + o(1).$$

Now setting $F := \{\Phi_{B_j} \leq b^+(\alpha, \beta, \eta, N)\}$, $G := \{x \in \mathcal{H}_N(\alpha)\}$ and using Chebyshev's inequality we obtain

$$\mathbf{P}(C \cap F|G) \leq \mathbf{P}(G)^{-1} N^{-2\beta_j(1-\alpha^2)-\varepsilon/2} \mathbb{E}(\mathbb{1}_F \mathbb{1}_G|\mathcal{H}_N(\alpha) \cap D_j|)$$

$$(4.2) \qquad = \mathbf{P}(G)^{-1} N^{-2\beta_j(1-\alpha^2)-\varepsilon/2} \mathbb{E}\left(\sum_{y \in D_j} \mathbb{1}_{\{x,y \in \mathcal{H}_N(\alpha)\}} \mathbb{1}_F\right)$$

$$\leq \mathbf{P}(G)^{-1} N^{2\beta_j\alpha^2-\varepsilon/2} \sup_{y \in D_j} \mathbf{P}(\{x, y \in \mathcal{H}_N(\alpha)\} \cap F).$$

Now by Lemma 2.1 and (2.15), we know that

$$(4.3) \qquad\qquad \mathbf{P}(G)^{-1} \leq N^{2\alpha^2+\varepsilon/8}$$

for $N$ large enough. Moreover, by Lemma 8.3, since $\gamma^* = 2/(2 - \beta_j) > 1$, for $\eta$ and $K$ small and large enough, respectively, we have

$$\sup_{y \in D_j} \mathbf{P}(\{x, y \in \mathcal{H}_N(\alpha)\} \cap F) \leq N^{-2\alpha^2 F_{2,\beta_j}(1)+\varepsilon/8}$$

$$= N^{-2\alpha^2(1+\beta_j)+\varepsilon/8}.$$

Plugging this together with (4.3) back in (4.2) yields

$$\mathbf{P}(C \cap F|G) \leq N^{2\alpha^2+\varepsilon/8+2\beta_j\alpha^2-\varepsilon/2-2\alpha^2(1+\beta_j)+\varepsilon/8}$$

$$= N^{-\varepsilon/4} \to 0$$

as $N$ goes to infinity. This concludes the proof of the upper bound. $\quad\square$

## 5. Theorem 1.6: Pairs of high points.

Let us first explain the idea behind the result (which is inspired by [5]). We partition $V_N$ into boxes of side length $N^\beta$, and let $\gamma > 0$. By a rerun of Theorem 1.3, the number of boxes $B$ in this collection such that

$$(5.1) \qquad\qquad \Phi_B \geq \gamma\alpha(1 - \beta)2\sqrt{g} \log N$$

is roughly

$$(5.2) \qquad \begin{cases} N^{m_\gamma}, & \text{if } m_\gamma := 2(1 - \beta)(1 - (\alpha\gamma)^2) > 0, \\ 0, & \text{otherwise.} \end{cases}$$

By conditioning over $\mathcal{F}_{\partial B}$ and using again Theorem 1.3, we obtain that each such box will approximately contain

$$(5.3) \qquad N^{4\beta(1-\eta_\gamma^2)}, \qquad \eta_\gamma := \alpha^2\left(\frac{1-\gamma(1-\beta)}{\beta}\right)^2$$

pairs of $\alpha$-high points. Multiplying (5.2) and (5.3) gives us a total of

$$(5.4) \qquad N^{2+2\beta-2\alpha^2 F_{2,\beta}(\gamma)}$$

pairs. Maximizing (5.4) over all $\gamma$ for which $m_\gamma > 0$ yields $N^{\rho(\alpha,\beta)}$ pairs of $\alpha$-high points, as claimed.

PROOF OF THE UPPER BOUND. We first rewrite (1.12). If we denote by $\gamma_m$ the value of $\gamma$ which achieves the minimum in (1.12), then by monotonicity it follows that $\gamma_m = \min\{\gamma_*, \gamma_+\}$, where $\gamma_*$ is the unconstrained minimizer of $F_{2,\beta}$ and where

$$\begin{aligned}
\gamma_+ &= \sup\{\Gamma_{\alpha,\beta}\} \\
&= \sup\{\gamma : 2 - 2\beta - 2\alpha^2\{\gamma^2(1-\beta)\} \geq 0\} \\
&= \frac{1}{\alpha} \geq 1.
\end{aligned}$$

Thus

$$\begin{aligned}
\rho(\alpha,\beta) &\geq 2 + 2\beta - 2\alpha^2 F_{2,\beta}(1) \\
&= 2 + 2\beta - 2\alpha^2((1-\beta) + 2\beta) \\
&= 2(1-\alpha^2)(1+\beta).
\end{aligned}$$

Now by Theorem 1.3, for $\lambda > 0$ the number of pairs of high points within distance $N^{\lambda\beta}$ of each other is at most $N^{2(1-\alpha^2)+2\beta\lambda}$. For this last quantity to be less than $N^{\rho(\alpha,\beta)}$ it suffices that

$$2(1-\alpha^2) + 2\beta\lambda \leq 2(1-\alpha^2)(1+\beta),$$

that is, that $\lambda \leq (1-\alpha^2)$. Moreover, one easily verifies that $\rho(\alpha,\beta)$ is increasing in $\beta$. Therefore to prove the upper bound we only need to show that for all $\delta > 0$ there exists $h$ such that for all $\beta' \in [\beta(1-\alpha^2), \beta]$

$$\mathbf{P}(|\{x,y\} \in \mathcal{H}_N(\alpha) : N^{\beta'h} \leq d(x,y) \leq N^{\beta'}| \geq N^{\rho(\alpha,\beta')+\delta}) \to 0$$

as $N$ tends to infinity. To prove this claim we will argue separately depending on whether or not $\gamma_*(\alpha,\beta') = \gamma_m(\alpha,\beta')$. Consider the former situation, that is, $\gamma_*(\alpha,\beta') = \gamma_m(\alpha,\beta')$. Let

$$E := \{|\{x,y\} \in \mathcal{H}_N(\alpha) : N^{\beta'h} \leq d(x,y) \leq N^{\beta'}| \geq N^{\rho(\alpha,\beta')+\delta}\}.$$

Then by Chebyshev's inequality

$$\mathbf{P}(E) \le N^{-\rho(\alpha,\beta')-\delta}\mathbb{E}\left(\sum_{x,y:\, N^{\beta'h}\le d(x,y)\le N^{\beta'}}\mathbb{1}_{\{x,y\in\mathcal{H}_N(\alpha)\}}\right)$$

$$\le N^{-\rho(\alpha,\beta')-\delta}N^2 N^{2\beta'} N^{-2\alpha^2 F_{2,\beta'}(\gamma^*)+\delta/2} = N^{-\delta/2}$$

where by Lemma 8.2 the last inequality holds provided $h$ is close enough to 1.

Let us now consider the case where $\gamma_* > \gamma_m$. For each box in $\Pi_{\beta'}$, we create a bigger box by adding the eight contiguous boxes. Let $\mathcal{B}$ be this collection (of boxes). For each of these boxes, associate a concentric box two times bigger. Let $\mathcal{C}$ be this new collection. The idea behind this construction is that any pair of two points within distance $N^{\beta'}$ is contained in at least one box in $\mathcal{B}$. Now let $\varepsilon > 0$ and

$$D := \left\{\max_{C\in\mathcal{C}}\Phi_C \le \alpha(\gamma_m+\varepsilon)(1-\beta')2\sqrt{g}\log N\right\}.$$

Then, because $\gamma_m = \gamma_+ = 1/\alpha$, $\mathbf{P}(D^c) \to 0$ [by the union bound, Lemma 2.1 and (2.6)]. Thus,

$$\mathbf{P}(E) = o(1) + \mathbf{P}(E\cap D)$$

$$\le o(1) + N^{-\rho(\alpha,\beta')-\delta}N^2 N^{2\beta'} N^{-2\alpha^2 F_{2,\beta'}(\gamma_m+\varepsilon)+\delta/2}$$

for $h$ sufficiently close to 1. In the last step we have used Chebyshev's inequality and Lemma 8.3. Since $2 + 2\beta' - 2\alpha^2 F_{2,\beta'}(\gamma_m+\varepsilon)$ tends to $\rho(\alpha,\beta')$ as $\varepsilon$ tends to 0, by choosing $\varepsilon$ small enough we see that $\mathbf{P}(E) \to 0$, which completes the proof of the upper bound. □

PROOF OF THE LOWER BOUND.    Let $\gamma < \gamma_+$, and let

$$C := \{|\{(x,y)\in\mathcal{H}_N(\alpha):d(x,y)\le N^\beta\}|\le N^{\rho(\alpha,\beta)-\delta}\}.$$

Let $m_\gamma := 2 - 2\beta - 2\alpha^2 F_{0,\beta}(\gamma)$,

$$F := \{B\in\Pi_\beta:\Phi_B \ge 2\gamma(1-\beta)\alpha\sqrt{g}\log N\}$$

and $D := \{\#F \ge N^{m_\gamma-\delta/2}\}$. Since $m_\gamma = 2(1-\beta)(1-\gamma^2\alpha^2)$, the proof of Theorem 1.3 shows that $\mathbf{P}(D^c) \to 0$. Hence

$$\mathbf{P}(C) = o(1) + \mathbf{P}(C\cap D).$$

On the event $D$, the set $F$ contains at least $N^{m_\gamma-\delta/2}$ boxes $\{B_j\}$. Let $D_j := \{\Phi_{B_j} \ge 2\gamma(1-\beta)\alpha\sqrt{g}\log N\}$. Notice that

$$C\cap D \subset E := \bigcup_{j=1}^{N^{m_\gamma-\delta/2}}\{|\mathcal{H}_N(\alpha)\cap B_j|\le N^{(\rho(\alpha,\beta)-m_\gamma)/2-\delta/8}\}\cap D_j.$$

Fixing $\eta > 0$ let

$$A := \bigcup_{B \in \Pi_\beta} \bigcup_{y \in B(x_B, N^\beta/2)} \{ |\mathbb{E}(\Phi_y | \mathcal{F}_\beta) - \Phi_B| \geq \eta 2\sqrt{g} \log N \}.$$

Then by [2], Lemma 12, $\mathbf{P}(A)$ tends to 0. So

$$\mathbf{P}(C \cap D)$$

$$\leq o(1) + \mathbf{P}(A^c \cap E)$$

$$\leq o(1)$$

$$+ N^{m_\gamma - \delta/2} \mathbf{P}(|\mathcal{H}_{N^\beta}^{1/4}((\alpha(1 - \gamma(1 - \beta)) + \eta)/\beta)| \leq N^{(\rho(\alpha,\beta) - m_\gamma)/2 - \delta/8}).$$

Since

$$(\rho(\alpha, \beta) - m_\gamma)/2 = 2\beta\left(1 - \frac{\alpha^2(1 - \gamma(1 - \beta))^2}{\beta^2}\right),$$

by Theorem 1.3 the probability in the right-hand side is lower than $e^{-c(\log N)^2}$ for some $c > 0$ as soon as $\eta$ is small enough. Therefore $\mathbf{P}(C \cap D)$ tends to 0, which completes the proof. $\square$

## 6. Theorem 1.7: The biggest high square. Our proof is inspired by [4].

PROOF OF THE LOWER BOUND. Here we use the notation introduced in the construction of Section 2. Let $1/2 < \alpha < 1$, and take $k$ such that

$$(6.1) \qquad \alpha_k := \alpha(K - k + 1)/K > 1/2 - \eta/2 - \delta.$$

Let $\underline{B}^{(k)}$ be any element of $D_k$ [which is defined in (2.8)], let $B_{1,k}$ be the last box of $\underline{B}^{(k)}$, and let $B$ be an $N^{\alpha_k}/2$ sized box centered in $B_{1,k}$. Let

$$A := \bigcup_{y \in B} \{ |\mathbb{E}(\Phi_y | \mathcal{F}_{\alpha_k}) - \mathbb{E}(\Phi_{x_B} | \mathcal{F}_{\alpha_k})| \geq \varepsilon 2\sqrt{g}(\alpha - \alpha_k)(1 - \gamma_k) \log N \},$$

where $\varepsilon > 0$. Since $\mathbf{P}(A^c)$ and $\mathbf{P}(C_k)$ tend to 1 (by [2], Lemma 12 and (2.10), (2.16), resp.],

$$\mathbf{P}(D_N(\eta) \leq N^{1/2 - \eta/2 - \delta})$$

$$\leq o(1) + \mathbf{P}\left( C_k \cap \left\{ \min_{x \in B} \Phi_x \leq \eta 2\sqrt{g} \log N \right\} \cap A^c \right)$$

and

$$\mathbf{P}\left( C_k \cap \left\{ \min_{x \in B} \Phi_x \leq \eta 2\sqrt{g} \log N \right\} \cap A^c \right)$$

$$\leq \mathbf{P}\left( \min_{x \in B} (\Phi_x - \mathbb{E}(\Phi_x | \mathcal{F}_{\alpha_k})) \leq 2\sqrt{g} \log N(\eta - (\alpha - \alpha_k)(1 - \gamma_K)(1 - \varepsilon)) \right)$$

$$\leq \mathbf{P}\left( \sup_{V_{N^{\alpha_k}}^{1/4}} \Phi \geq 2\sqrt{g} \log N((\alpha - \alpha_k)(1 - \gamma_K)(1 - \varepsilon) - \eta) \right).$$

If

$$2\sqrt{g}\log N((\alpha-\alpha_k)(1-\gamma_K)(1-\varepsilon)-\eta) > 2\sqrt{g}\log N^{\alpha_k}$$

for $N$ large enough, then [2], Theorem 2a completes the proof. But this last condition, together with (6.1), holds provided that we choose $1/2 + \eta/2 < k/K < 1/2 + \eta/2 + \delta$ and $\alpha$ close to 1 as well as $\varepsilon$ small enough and $K, N$ large enough. $\square$

PROOF OF THE UPPER BOUND. Let $\alpha = 1/2 - \eta/2$ and fix $\beta = \alpha + \delta$. First remark that

$$(6.2) \qquad \mathbf{P}\left(\bigcup_{B\in\Pi_\beta}\{\Phi_B \geq 2\sqrt{g}(1-\alpha)\log N\}\right) \xrightarrow{N\to\infty} 0$$

follows easily from the union bound, Lemma 2.1 and (2.6). Next let

$$F := \left\{\bigcap_{B\in\Pi_\beta}\{\Phi_B \leq 2\sqrt{g}(1-\alpha)\log N\}\right\},$$

$$C := \left\{\bigcup_{B\in\Pi_\beta}\{\forall\, x\in B,\ \Phi_x \geq \eta 2\sqrt{g}\log N\}\right\}.$$

Then for $N$ large enough

$$(6.3) \qquad \begin{aligned} \mathbf{P}(D_N(\eta) \geq N^{\alpha+2\delta}) &\leq \mathbf{P}(C) \\ &\leq \mathbf{P}(F^c) + \mathbf{P}(C \cap F) \\ &= o(1) + \mathbb{E}(\mathbf{P}(C|\mathcal{F}_\beta)\mathbb{1}_F) \end{aligned}$$

where we have used (6.2). For any $B \in \Pi_\beta$ we let $B^{(1/4)}$ be the sub-box $B(x_B, N^\beta/2)$, and for $\varepsilon > 0$ we define

$$A := \bigcup_{B\in\Pi_\beta}\ \bigcup_{y\in B^{(1/4)}}\{|\mathbb{E}(\Phi_y|\mathcal{F}_{\partial B}) - \Phi_B| \geq \varepsilon 2\sqrt{g}\log N\}.$$

Then by [2], Lemma 12, $\mathbf{P}(A)$ tends to 0. Thus (6.3) becomes

$$(6.4) \qquad \mathbf{P}(D_N(\eta) \geq N^{\alpha+2\delta}) \leq o(1) + \mathbb{E}(\mathbf{P}(C|\mathcal{F}_\beta)\mathbb{1}_F\mathbb{1}_{A^c}).$$

Now remark that

$$\mathbf{P}(C|\mathcal{F}_\beta) \leq \left(\frac{N}{N^\beta}\right)^2\max_{B\in\Pi_\beta}\mathbf{P}(\forall\, x\in B,\ \Phi_x \geq \eta 2\sqrt{g}\log N|\mathcal{F}_\beta).$$

But on $A^c \cap F$,

$$\begin{aligned} &\mathbf{P}(\forall\, x\in B,\ \Phi_x \geq \eta 2\sqrt{g}\log N|\mathcal{F}_\beta) \\ &\qquad \leq \mathbf{P}(\forall\, x\in B,\ (\Phi_x - \mathbb{E}(\Phi_x|\mathcal{F}_\beta)) \geq (\eta - (1-\alpha+\varepsilon))2\sqrt{g}\log N|\mathcal{F}_\beta) \\ &\qquad = \mathbf{P}\left(\max_{V_{N^\beta}^{1/4}}\Phi \leq 2\sqrt{g}(\alpha+\varepsilon)\log N\right). \end{aligned}$$

By [2], Theorem 2, this last quantity can be bounded by $\exp(-d(\log N)^2)$ for some $d > 0$ provided $\beta > \alpha + \varepsilon$, that is, $\varepsilon < \delta$. Consequently, on $A^c \cap F$,

$$\mathbf{P}(C|\mathcal{F}_\beta) \leq \exp[2(1-\beta)\log N - d(\log N)^2]$$

which tends to 0 as $N \to \infty$. In view of (6.4) this concludes the proof of the upper bound. $\quad\square$

## 7. Results on the conditioned free field.

### 7.1. *Width of downward spike.*

PROOF OF THEOREM 1.1. We start by establishing the lower bound. Fix $\delta, \varepsilon > 0$ and let $a := 2\sqrt{g} + \varepsilon$. By FKG (e.g., see [8], (B.13)) we have

$$\mathbf{P}(\cdot|\Omega^+_{N,l}) \prec \mathbf{P}^a(\cdot|\Omega^+_{N,l})$$

where $\mathbf{P}^a$ denotes the Gaussian measure with covariance $G_N$ and mean $a \log N$. So

$$
\begin{aligned}
(7.1) \quad &\mathbf{P}\left(\frac{\log D^+_N(\eta)}{\log N} \leq \frac{1}{2} - \frac{\eta}{2} - \delta \Big| \Omega^+_{N,l}\right) \\
&\leq \mathbf{P}^a\left(\frac{\log D^+_N(\eta)}{\log N} \leq \frac{1}{2} - \frac{\eta}{2} - \delta \Big| \Omega^+_{N,l}\right) \\
&\leq \frac{\mathbf{P}^a(\log D^+_N(\eta)/\log N \leq 1/2 - \eta/2 - \delta)}{\mathbf{P}^a(\Omega^+_{N,l})}.
\end{aligned}
$$

Now by [2], Theorem 2, the denominator of (7.1) tends to 1 while by Theorem 1.7 the numerator tends to 0 provided $\varepsilon$ is small enough.

Turning to the upper bound, we fix $\delta > 0$ and let $\alpha = 1/2 - \eta/2$, $\mathcal{B} := \Pi_{\alpha+\delta}$, $\mathbf{P}^+(\cdot) := \mathbf{P}(\cdot|\Omega^+_{N,l})$ and $\mathbb{E}^+(\cdot) := \mathbb{E}(\cdot|\Omega^+_{N,l})$. For $\varepsilon > 0$, $x \in V^l_N$,

$$
\begin{aligned}
(7.2) \quad &\mathbb{E}^+(\Phi_x) \geq (2\sqrt{g}(1-\varepsilon)\log N)\mathbf{P}^+(\Phi_x \geq 2\sqrt{g}(1-\varepsilon)\log N) \\
&\geq (2\sqrt{g}(1-\varepsilon)\log N) \\
&\quad \times \left(1 - \sup_{y \in V^l_N} \mathbf{P}^+(|\Phi_y - 2\sqrt{g}\log N| \geq \varepsilon 2\sqrt{g}\log N)\right).
\end{aligned}
$$

By [2], Theorem 4, the supremum in (7.2) tends to 0 as $N$ tends to infinity, and therefore $\inf_{y \in V^l_N} \mathbb{E}^+(\Phi_y)$ is greater than $(2\sqrt{g} - \varepsilon)\log N$ for any $\varepsilon$ provided $N$ is large enough. For any $B \in \mathcal{B}$, $\Phi_B$ is a convex linear combination of such $\Phi_x$ (with positive coefficients), and consequently we also have

$$(7.3) \qquad \inf_{B \in \mathcal{B}} \mathbb{E}^+(\Phi_B) \geq 2\sqrt{g}(1-\varepsilon)\log N$$

as soon as $N$ is large enough. Next, using the Brascamp–Lieb inequality (see [8], (B.21)) in the first step we obtain that

$$\mathbf{P}^+(\Phi_B - \mathbb{E}^+(\Phi_B) \leq -(1-\alpha-\varepsilon)2\sqrt{g}\log N)$$

(7.4)
$$\leq e^{-(1-\alpha-\varepsilon)^2 4g(\log N)^2/(2\mathrm{Var}\Phi_B)}$$

$$\leq e^{-(2(1-\alpha-\varepsilon)^2/(1-\alpha-\delta/2))\log N}$$

$$\leq e^{-2(1-\alpha+\nu)\log N}$$

for some $\nu > 0$ provided $\varepsilon$ is small enough. On the other hand, if we let

$$C := \bigcup_{B \in \mathcal{B}} \{\Phi_B \leq 2\sqrt{g}\alpha\log N\},$$

then by (7.3)

$$C \subset \bigcup_{B \in \mathcal{B}} \{\Phi_B - \mathbb{E}^+(\Phi_B) \leq 2\sqrt{g}(\alpha-1+\varepsilon)\log N\}$$

and therefore by the union bound and (7.4), $\lim_{N\to\infty}\mathbf{P}^+(C) = 0$. Next, if for every box $B \in \mathcal{B}$ we let $B^{(\varepsilon')} := B(x_B, \varepsilon' N^{\alpha+\delta})$ and define

$$F := \bigcup_{B \in \mathcal{B}} \bigcup_{x \in B^{(\varepsilon')}} \{|\Phi_B - \mathbb{E}(\Phi_x|\mathcal{F}_{\partial B})| \geq 2\sqrt{g}\varepsilon''\log N\},$$

then we know from [2], Lemma 12, and [2], Theorem 3 that $\mathbf{P}^+(F) \to 0$ for any choice of $\varepsilon'' > 0$, as long as $\varepsilon' > 0$ is small enough. Now for $N$ big enough, since $1-\eta = 2\alpha$,

$$\mathbf{P}^+(\log D_N^+(\eta) \geq (\alpha+2\delta)\log N)$$

$$\leq \mathbf{P}^+\left(\bigcup_{B \in \mathcal{B}} \sup_{x \in B^{(\varepsilon')}} \Phi_x \leq 2\sqrt{g}(2\alpha)\log N\right)$$

$$= o(1) + \mathbf{P}^+\left(C^c \cap F^c \cap \bigcup_{B \in \mathcal{B}} \sup_{x \in B^{(\varepsilon')}} \Phi_x \leq 2\sqrt{g}(2\alpha)\log N\right)$$

$$\leq o(1) + \mathbf{P}^+\left(\bigcup_{B \in \mathcal{B}} \sup_{x \in B^{(\varepsilon')}} (\Phi_x - \mathbb{E}(\Phi_x|\mathcal{F}_{\partial B})) \leq 2\sqrt{g}(\alpha+\varepsilon'')\log N\right)$$

$$\leq o(1) + N^{2-2\alpha} \max_{B \in \mathcal{B}} \mathbf{P}^+\left(\sup_{x \in B^{(\varepsilon')}} (\Phi_x - \mathbb{E}(\Phi_x|\mathcal{F}_{\partial B}))\right.$$
$$\left. \leq 2\sqrt{g}(\alpha+\varepsilon'')\log N\right)$$

$$= o(1) + N^{2-2\alpha} \max_{B \in \mathcal{B}} \mathbb{E}^+\left(\mathbf{P}^+\left(\sup_{x \in B^{(\varepsilon')}} (\Phi_x - \mathbb{E}(\Phi_x|\mathcal{F}_{\partial B}))\right.\right.$$
$$\left.\left. \leq 2\sqrt{g}(\alpha+\varepsilon'')\log N \Big| \mathcal{F}_{B^c}\right)\right).$$

Now $\mathbf{P}^+(\cdot|\mathcal{F}_{B^c})$ is equal to $\mathbf{P}(\cdot|\mathcal{F}_{B^c})$ conditioned on being positive (inside $B$). Thus, $\mathbf{P}^+(\cdot|\mathcal{F}_{B^c}) \succ \mathbf{P}(\cdot|\mathcal{F}_{B^c})$ (see, e.g., [8], (B.5)). Therefore the right-hand side of the last equation is less than

$$o(1) + N^{2-2\alpha} \max_{B \in \mathcal{B}} \mathbb{E}^+ \bigg( \mathbf{P} \bigg( \sup_{x \in B^{(\varepsilon')}} (\Phi_x - \mathbb{E}(\Phi_x|\mathcal{F}_{\partial B}))$$
$$\leq 2\sqrt{g}(\alpha + \varepsilon'') \log N \Big| \mathcal{F}_{B^c} \bigg) \bigg)$$

which is itself less than $e^{-c(\log N)^2}$ for some $c > 0$ by [2], Theorem 2, provided $\varepsilon'' < \delta$. This completes the proof of the lower bound. $\square$

### 7.2. *Number of low points.*

PROOF OF PROPOSITION 1.2. As for Theorem 1.1, the lower bound follows easily from FKG, combined with Theorem 1.3 in this case. To prove the upper bound, remark again that by [2], Theorem 4,

$$\inf_{x \in V_N^l} \mathbb{E}^+(\Phi_x) \geq (2\sqrt{g} - \varepsilon) \log N$$

for any $\varepsilon > 0$, provided $N$ is large enough. Combined with the Brascamp–Lieb inequality, this implies that for some $C > 0$

$$\mathbf{P}^+(\Phi_x \leq 2\sqrt{g}(1-\eta) \log N) \leq C \exp[-2(\eta - \varepsilon)^2 \log N]$$

uniformly on $x \in V_N^l$ provided $N$ is large enough. An application of the union bound completes the proof. $\square$

## 8. Gaussian computations.

When $x$ is an $\alpha$-high point, the next lemma gives the typical value of the GFF for points within distance $N^\beta$ of $x$.

LEMMA 8.1. *Let $B := B(x, N^\beta)$, $\varepsilon > 0$, $b^\pm(\alpha, \beta, \varepsilon, N) = 2\sqrt{g}(\alpha(1-\beta) \pm \varepsilon) \log N$ and $I(\alpha, \beta, \varepsilon, N) := [b^-(\alpha, \beta, \varepsilon, N), b^+(\alpha, \beta, \varepsilon, N)]$. Then*

$$\tag{8.1} \max_{x \in V_N^l} \mathbf{P}(\Phi_x \notin I(\alpha, \beta, \varepsilon, N) | \Phi_x \geq 2\sqrt{g}\alpha \log N) \longrightarrow 0$$

*as $N \to \infty$.*

PROOF. To simplify notation we write $I$, $b^+$ and $b^-$ instead of $I(\alpha, \beta, \varepsilon, N)$, $b^+(\alpha, \beta, \varepsilon, N)$ and $b^-(\alpha, \beta, \varepsilon, N)$. For $\eta > 0$ we have [by (2.6) and (2.15)]

$$\tag{8.2} \mathbf{P}(\Phi_x \geq 2\sqrt{g}\alpha(1+\eta) \log N | \Phi_x \geq 2\sqrt{g}\alpha \log N) \xrightarrow{N \to \infty} 0.$$

Thus

$$\mathbf{P}(\Phi_B \notin I | \Phi_x \geq 2\sqrt{g}\alpha \log N)$$
$$= o(1) + \mathbf{P}(\Phi_B \notin I, \Phi_x \leq 2\sqrt{g}\alpha(1+\eta)\log N | \Phi_x \geq 2\sqrt{g}\alpha \log N)$$
$$\leq o(1) + \mathbf{P}(\Phi_B \notin I | \Phi_x \in (1, 1+\eta)2\sqrt{g}\alpha \log N).$$

Thus it suffices to prove that

$$\mathbf{P}(\Phi_B \notin I | \Phi_x \in (1, 1+\eta)2\sqrt{g}\alpha \log N)$$

tends to 0. Now since $\Phi_x = \Phi_x - \Phi_B + \Phi_B$ and $\Phi_x - \Phi_B$ is independent of $\Phi_B$ we have $\mathrm{Cov}(\Phi_x, \Phi_B) = \mathrm{Var}(\Phi_B)$. Hence

$$\Phi_B = \frac{\mathrm{Var}(\Phi_B)}{\mathrm{Var}(\Phi_x)}\Phi_x + Z$$

where $Z$ is centered Gaussian and independent of $\Phi_x$. Now

$$\frac{\mathrm{Var}(\Phi_B)}{\mathrm{Var}(\Phi_x)} = (1-\beta) + O\left(\frac{1}{\log N}\right)$$

and therefore $\mathrm{Var}(Z) = O(\log N)$. Consequently

$$\mathbf{P}(\Phi_B \geq b^+ | \Phi_x \in (1, 1+\eta)2\sqrt{g}\alpha \log N)$$
$$\leq \mathbf{P}\left(Z + \left((1-\beta) + O\left(\frac{1}{\log N}\right)\right)\alpha(1+\eta)2\sqrt{g}\log N \geq b^+\right)$$

tends to 0 when $\eta$ is small enough. Similarly

$$\mathbf{P}(\Phi_B \leq b^- | \Phi_x \in (1, 1+\eta)2\sqrt{g}\alpha \log N)$$
$$\leq \mathbf{P}\left(Z + \left((1-\beta) + O\left(\frac{1}{\log N}\right)\right)\alpha 2\sqrt{g}\log N \leq b^-\right) \to 0.$$

This completes the proof of the lemma. $\square$

LEMMA 8.2. *Let* $0 < \alpha < \beta < 1$, $\delta > 0$ *and*

$$S := \{(x, y) \in V_N^l : N^{\beta(1-\varepsilon)} \leq d(x, y) \leq N^\beta\}.$$

*Then there exists* $C, \varepsilon_0 > 0$ *such that for* $\varepsilon \leq \varepsilon_0$ *and all* $N$

$$\max_{(x,y)\in S} \mathbf{P}(x, y \in \mathcal{H}_N(\alpha)) \leq C N^{-2\alpha^2 F_{2,\beta}(\gamma^*) + \delta}$$

*where* $\gamma^* = 2/(2-\beta)$. *Moreover,* $\varepsilon_0$ *can be chosen uniformly in* $(\alpha, \beta)$ *over compact sets of* $(0, 1)^2$.

PROOF. Let $Z := \Phi_x + \Phi_y$ and notice that

$$\{x, y \in \mathcal{H}_N(\alpha)\} \subset \{Z \geq 4\alpha\sqrt{g}\log N\}.$$

We have $\mathrm{Var}(Z) = \mathrm{Var}(\Phi_x) + \mathrm{Var}(\Phi_y) + 2\,\mathrm{Cov}(\Phi_x, \Phi_y)$, and by Lemma 2.1

$$(8.3) \qquad \mathrm{Var}(\Phi_x) \leq g\log N + O(1),$$

$$(8.4) \qquad \mathrm{Var}(\Phi_y) \leq g\log N + O(1),$$

$$(8.5) \qquad \mathrm{Cov}(\Phi_x, \Phi_y) \leq g\log N - g\beta(1-\varepsilon)\log N + O(1).$$

Thus, $\mathrm{Var}(Z) \leq (2g(2-\beta) + O(\varepsilon) + O(1/\log N))\log N$. Therefore, since $F_{2,\beta}(\gamma^*) = \gamma^*$,

$$\mathbf{P}(Z \geq 4\alpha\sqrt{g}\log N)$$

$$\leq \exp\Big(-\frac{(4\alpha\sqrt{g}\log N)^2}{2(2g(2-\beta) + O(\varepsilon) + O(1/\log N))\log N}\Big)$$

$$= \exp(-2\alpha^2\gamma^*(1 + O(\varepsilon) + O(1/\log N))\log N)$$

$$\leq CN^{-2\alpha^2 F_{2,\beta}(\gamma^*) + O(\varepsilon)}$$

which concludes the proof. $\square$

LEMMA 8.3. *Let $0 < \alpha < \beta < 1$ and $\delta > 0$. For $(x, y) \in S$ (as defined in Lemma 8.2) we let $T(x, y)$ be the set of boxes $B$ of side length $2N^\beta$ whose centered sub-box of size $N^\beta$ contains $x$ and $y$. Then there exists $C, \varepsilon_0 > 0$ such that for $\varepsilon \leq \varepsilon_0$ and all $N$*

$$\max_{\substack{(x,y) \in S, \\ B \in T(x,y)}} \mathbf{P}(\{x, y \in \mathcal{H}_N(\alpha)\} \cap \{\Phi_B \leq \gamma\alpha(1-\beta)2\sqrt{g}\log N\})$$

$$\leq CN^{-2\alpha^2 F_{2,\beta}(\min\{\gamma^*, \gamma\}) + \delta}.$$

*Moreover, $\varepsilon_0$ can be chosen uniformly in $(\alpha, \beta)$ over compact sets of $(0,1)^2$.*

PROOF. We let

$$E := \{x, y \in \mathcal{H}_N(\alpha)\} \cap \{\Phi_B \leq \gamma\alpha(1-\beta)2\sqrt{g}\log N\},$$

and argue separately depending on whether or not $\gamma \geq \gamma^*$. When $\gamma \geq \gamma^*$, it suffices to notice that $\mathbf{P}(E) \leq \mathbf{P}(x, y \in \mathcal{H}_N(\alpha))$ and to conclude using Lemma 8.2 (since in this case $\min\{\gamma^*, \gamma\} = \gamma^*$). In the case $\gamma < \gamma^*$, we introduce $a := 1 - \gamma(1-\beta)$ and $b := \gamma(2-\beta) - 2$. Then since

$$\gamma < \gamma^* = \frac{2}{2-\beta} < \frac{1}{1-\beta},$$

we have $a > 0$ and $b < 0$. Consequently, if we let $Z := a(\Phi_x + \Phi_y) + b\Phi_B$, then

$$(8.6) \qquad E \subset \{Z \geq (2a + b\gamma(1-\beta))2\alpha\sqrt{g}\log N\}.$$

Now

$$(8.7) \qquad \begin{aligned} \mathrm{Var}(Z) = {} & a^2\mathrm{Var}(\Phi_x) + a^2\mathrm{Var}(\Phi_y) + b^2\mathrm{Var}(\Phi_B) \\ & + 2ab\,\mathrm{Cov}(\Phi_x, \Phi_B) + 2ab\,\mathrm{Cov}(\Phi_y, \Phi_B) + 2a^2\,\mathrm{Cov}(\Phi_x, \Phi_y). \end{aligned}$$

Upper bounds for $\mathrm{Var}(\Phi_x)$, $\mathrm{Var}(\Phi_y)$ and $\mathrm{Cov}(\Phi_x, \Phi_y)$ are given by (8.3), (8.4) and (8.5), respectively. By Lemma 2.1 we have

$$\mathrm{Var}(\Phi_B) = \mathrm{Var}(\Phi_{x_B}) - \mathrm{Var}_{\mathcal{F}_{\partial B}}(\Phi_{x_B}) \leq g(1-\beta)\log N + O(1).$$

As for the covariances,

$$\begin{aligned} \mathrm{Cov}(\Phi_x, \Phi_B) &= \mathbb{E}(\mathbb{E}(\Phi_x|\mathcal{F}_{\partial B})\mathbb{E}(\Phi_{x_B}|\mathcal{F}_{\partial B})) \\ &= \mathrm{Cov}(\Phi_x, \Phi_{x_B}) - \mathrm{Cov}_{\mathcal{F}_{\partial B}}(\Phi_x, \Phi_{x_B}) \\ &\geq g(\log N - \log|x - x_B|) - g(\beta\log N - \log|x - x_B|) + O(1) \\ &= g(1-\beta)\log N + O(1). \end{aligned}$$

Similarly, we have

$$\mathrm{Cov}(\Phi_y, \Phi_B) \geq g(1-\beta)\log N + O(1).$$

Plugging these estimates in (8.7) and setting

$$f(a, b, \beta) := 2a^2(2-\beta) + b^2(1-\beta) + 4(1-\beta)ab$$

yields

$$\mathrm{Var}(Z) \leq (f(a, b, \beta) + O(\varepsilon) + O((\log N)^{-1}))g\log N.$$

Now remark that $2a + b = \beta\gamma$ and therefore

$$4a^2 + b^2 + 4ab = (2a + b)^2 = \beta^2\gamma^2.$$

Thus

$$\begin{aligned} f(a, b, \beta) &= 4a^2 + b^2 + 4ab - \beta(2a^2 + b^2 + 4ab) \\ &= (4a^2 + b^2 + 4ab)(1-\beta) + 2\beta a^2 \\ &= (\beta\gamma)^2(1-\beta) + 2\beta a^2 \\ &= \beta(\beta\gamma^2(1-\beta) + 2a^2) \\ &= \beta((2a+b)(1-a) + 2a^2) \\ &= \beta(2a + b - ab), \end{aligned}$$

hence

$$(8.8) \qquad \mathrm{Var}(Z) \le (\beta(2a + b - ab) + O(\varepsilon) + O((\log N)^{-1}))g \log N.$$

Therefore, since $2a + b\gamma(1 - \beta) = 2a + b(1 - a) = 2a + b - ab$, (8.8) together with (8.6) yields

$$(8.9) \qquad \mathbf{P}(E) \le C \exp\left(\left(-2\alpha^2 \frac{2a + b - ab}{\beta} + O(\varepsilon)\right) \log N\right).$$

But remark that

$$\beta F_{2,\beta}(\gamma) = \beta\gamma^2(1 - \beta) + 2(1 - \gamma(1 - \beta))^2 = \beta\gamma^2(1 - \beta) + 2a^2$$
$$= (2a + b)(1 - a) + 2a^2 = 2a + b - ab.$$

In view of (8.9), this completes the proof. $\quad\square$

## 9. Further remarks and an open problem.

(i) In the same spirit as Theorem 1.1, we can show that Theorem 1.4 and Theorem 1.6 still hold in the case of the conditioned free field shifted (downward) by $2\sqrt{g} \log N$. The proof of these extensions is similar in spirit to the proof of Theorem 1.1. The upper bound follows from FKG, and the lower bound from adapting the proof of the original theorems.

(ii) It is easy to prove that Theorem 1.3 and Theorem 1.6 remain valid if one only counts the high points in $DN \cap \mathbb{Z}^2$ where $D$ is any open subset of $V := (0,1)^2$ which has positive distance to the complement of $V$. Only the lower bounds need to be adapted; one simply does so by restricting the families $\Pi$. to boxes in $DN$ instead of $V_N^l$.

(iii) Our results reveal that the extrema of the free field and those of its conditioned version (minus its mean) exhibit a similar behavior. It would therefore be interesting to further push the comparison between the two objects. The discrete two-dimensional free field can be seen as a random distribution, that is, a Gaussian sequence indexed by the family $C_0(V)$ of continuous real-valued functions with compact support in $V$. Indeed, $\Phi$ (which we momentarily denote $\Phi^N$) acts on such a function $f$ by

$$(\Phi^N, f) = \frac{1}{N^2} \sum_{x \in V_N} f(x/N)\Phi_x^N$$

and consequently for any $f$ and $g$ in $C_0(V)$

$$\mathrm{Cov}((\Phi^N, f), (\Phi^N, g)) = \frac{1}{N^4} \sum_{x,y \in V_N} f(x/N)g(y/N)G_N(x, y).$$

It can be proved using the invariance principle (cf. [1], Lemma 2.10) that the sequence $\{(\Phi^N, f), f \in C_0([0,1]^2)\}$ converges in law to the *continuous free* *field*, a continuous analogue of $\Phi^N$ which is presented in [12], for example. The continuous Gaussian free field can be simply defined as a (centered) Gaussian family $\{(\Phi, f)\}$ indexed by $\{f \in C_0(V)\}$ such that

$$\mathrm{Cov}((\Phi, f), (\Phi, g)) = \int_V f(x) g(y) \mathcal{G}(x, y) \, dx \, dy,$$

where $\mathcal{G}$ is the Green function of the Brownian motion killed as it exits $V$. Let $T_N : V_N^{\mathbb{R}} \to V_N^{\mathbb{R}}$ be defined by $(T_N(\Phi))_x = \Phi_x - \mathbb{E}^+(\Phi_x)$. We suspect that under $\mathbf{P}_N^+ o T_N^{-1}$, $\Phi_N$ might also converge weakly to the continuous Gaussian free field.

**Acknowledgments.** I am very grateful to Amir Dembo for bringing my attention to [2] in connection with his work on the two-dimensional simple random walk, and to Amir Dembo and Jean-Dominique Deuschel for suggesting Theorem 1.1 and Theorem 1.7.

DEPARTMENT OF MATHEMATICS
STANFORD UNIVERSITY
STANFORD, CALIFORNIA 94305
USA
E-MAIL: olivier.daviaud@gmail.com